\documentclass[10 pt, conference]{ieeeconf}
\IEEEoverridecommandlockouts
\usepackage{cite}
\usepackage{amsmath,amssymb,amsfonts,amsthm}
\usepackage{mathrsfs}
\usepackage{multirow}
\usepackage{algorithmic}
\usepackage{graphicx}
\usepackage{subcaption}
\usepackage{textcomp}
\usepackage{xcolor}
\usepackage{bm}
\usepackage[utf8]{inputenc}
\usepackage[T1]{fontenc}
\usepackage{booktabs}
\usepackage{makecell}
\def\BibTeX{{\rm B\kern-.05em{\sc i\kern-.025em b}\kern-.08em
		T\kern-.1667em\lower.7ex\hbox{E}\kern-.125emX}}

\usepackage{hyperref}

\usepackage{tikz}
\usepackage{pgfplots}
\usetikzlibrary{shapes.geometric, arrows, positioning, calc}

\tikzstyle{box} = [rectangle, minimum width=2.4cm, minimum height=1cm,text centered, draw=black]
\tikzstyle{arrow} = [thick,->]

\newcommand{\Rb}{{\mathbb{R}}}

\newcommand{\Ib}{{\mathbb{I}}}
\newcommand{\Zb}{{\mathbb{Z}}}

\newcommand{\Sb}{{\mathbb{S}}}

\newcommand{\bX}{{\bm{X}}}

\newcommand{\bc}{{\bm{c}}}
\newcommand{\bh}{{\bm{h}}}
\newcommand{\br}{{\bm{r}}}
\newcommand{\bg}{{\bm{g}}}
\newcommand{\bL}{{\bm{L}}}
\newcommand{\bM}{{\bm{M}}}

\newcommand{\Dc}{{\mathcal{D}}}

\newcommand{\Tc}{{\mathcal{T}}}

\newcommand{\Mc}{{\mathcal{M}}}

\newcommand{\Lc}{{\mathcal{L}}}

\newcommand{\tbc}{\tilde{\bm{c}}}

\newcommand{\tP}{ \tilde{P} }

\DeclareMathOperator{\tr}{{\rm tr}}

\DeclareMathOperator{\rank}{{rank}}
\DeclareMathOperator{\SOC}{{SOC}}
\DeclareMathOperator{\proj}{{proj}}

\DeclareMathOperator{\size}{{size}}

\theoremstyle{plain}
\newtheorem{problem}{Problem}
\newtheorem{theorem}{Theorem}

\newtheorem{remark}{Remark}

\begin{document}

	\title{Linear Model Predictive Control under Continuous Path Constraints via Parallelized Primal-Dual Hybrid Gradient Algorithm\\
		\thanks{$^{1}$Zishuo Li, Bo Yang, Jiayun Li, and Yilin Mo are with the Department of Automation, Tsinghua University, Beijing, 100084, China. 
			{\small \{lizs19,yang-b21,lijiayun22\}@mails.tsinghua.edu.cn,
				ylmo@tsinghua.edu.cn}}
		\thanks{$^{2}$Jiaqi Yan is with Department of Computer Science, Tokyo Institute of Technology, Tokyo, Japan. 
			{\small jyan@sc.dis.titech.ac.jp}
		}
	}
	\author{Zishuo Li$^{1}$, Bo Yang$^{1}$, Jiayun Li$^{1}$, Jiaqi Yan$^{2}$, Yilin Mo$^{1}$}
	
	\maketitle
	
	\begin{abstract}
		In this paper, we consider a Model Predictive Control~(MPC) problem of a continuous-time linear time-invariant system subject to continuous-time path constraints on the states and the inputs. By leveraging the concept of differential flatness, we can replace the differential equations governing the system with linear mapping between the states, inputs, and flat outputs~(including their derivatives). The flat outputs are then parameterized by piecewise polynomials, and the model predictive control problem can be equivalently transformed into a Semi-Definite Programming~(SDP) problem via Sum-of-Squares (SOS), ensuring constraint satisfaction at every continuous-time interval. We further note that the SDP problem contains a large number of small-size semi-definite matrices as optimization variables. To address this, we develop a Primal-Dual Hybrid Gradient~(PDHG) algorithm that can be efficiently parallelized to speed up the optimization procedure. Simulation results on a quadruple-tank process demonstrate that our formulation can guarantee strict constraint satisfaction, while the standard MPC controller based on the discretized system may violate the constraint inside a sampling period. Moreover, the computational speed superiority of our proposed algorithm is collaborated by numerical simulation.
	\end{abstract}
	
	\title{Linear Model Predictive Control under Continuous Path Constraints via Parallelized Primal-Dual Hybrid Gradient Algorithm\\
		\thanks{Identify applicable funding agency here. If none, delete this.}
	}
	\author{Zishuo Li, Bo Yang, Jiayun Li, Jiaqi Yan, Yilin Mo}
	
	
	\section{Introduction}

The optimal control theory aims to find control laws for a dynamical system in order to optimize a given objective function, which finds numerous applications in fields of engineering~\cite{ben2010optimal,sharp2011vehicle,satici2013robust} and economics~\cite{weber2011optimal,aseev2012infinite,kellett2019feedback} etc. Closed-form optimal control law can be found for certain unconstrained problems, such as linear-quadratic control problem~\cite{bertsekas2012dynamic}, or brachistochrone problem~\cite{clarke2013functional}. However, analytically solving the optimal control problem of continuous-time systems remains a challenging task. Furthermore, a vast majority of real-world dynamical systems operate under various constraints, such as input saturation or safety constraint on the state. For constrained optimal control problem, Pontryagin's maximum principle~\cite{hartl1995survey} can be used to derive necessary condition for optimality. However, in practice, only a small number of problems can be solved analytically. Therefore, algorithms, such as model predictive control, discretize the system and thus reducing the search space of the control input from the infinite dimensional function space into a finite dimensional space, where numerical optimization can be used. 

Dynamic Matrix Control~(DMC)~\cite{lundstrom1995limitations} and Model Algorithmic Control~(MAC)~\cite{rouhani1982model} are two formulations of MPC algorithm for discretized optimal control problems with constraints~\cite{garcia1989model}. Both formulations employ a zero-order hold for the control inputs, which implies that the control inputs are step functions and hence reside in a finite dimensional space. 
However, the discretization of a continuous time system means that one can only guarantee constraint satisfaction at all discrete-time instant, where constraint violation can occur in between. 

For control applications with high safety requirements, constraints violations can be intolerable. In order to meet the constraints at all time, Semi-Infinite Programming~(SIP)~\cite{djelassi2021recent} has been used to deal with infinite number of constraints. Several approaches for solving SIP have been proposed, and a common framework is to check constraint violations in intervals, and adaptively add additional discrete-time points until the tolerance level is guaranteed or no constraint violations occur. Chen et al~\cite{chen2005inequality} introduce $\epsilon$-tolerance on inequality constraints, which means that the constraints may still be violated up to $\epsilon$. Fu et al.~\cite{fu2015local} tighten the inequality constraints at discrete-time instant, hence guarantee the satisfactory of constraints over the whole interval. However, tighter constraints may lead to relatively conservative solution.

To address these issues, we parameterize the flat output of the continuous-time linear system by piecewise polynomials. The differential equation of the dynamic system is eliminated and replaced by flatness map between flat output $y$ and system state $x$, input $u$~\cite{763209}. In this way, the decision variables become finite-dimensional polynomial coefficients. On the other hand, the inequality constraints become non-negative polynomials over intervals, which are still infinite-dimensional. Fortunately, we can leverage Markov-Luk\'{a}cs theorem~\cite{Roh_SOS,szego_orth_poly} to transcribe a polynomial inequality constraint on an interval into an equivalent matrix Positive Semi-Definite (PSD) constraint, thus ensures the path constraints hold at every time interval. With this procedure, the continuous-time MPC problem can be transcribed into a Semi-Definite Programming (SDP) problem.

It is worth noticing that the SDP problem we formulate contains a large amount of small symmetric matrices. As a result, we propose to use parallel computing to speed up the calculation. To this end, we use a customized Primal-Dual Hybrid Gradient (PDHG) algorithm to solve the SDP problem. PDHG, also known as Chambolle-Pock~\cite{Chambolle2011}, is a well-known first-order algorithm dealing with convex optimization problems with equality constraints. 
For large scale problems, it has been one of the preferred first-order algorithm~\cite{NEURIPS2021_a8fbbd3b} due to the fact that it can be easily parallelized.
	
	The main contributions of this article are as follows.
	\begin{itemize}
		\item An equivalent formulation of continuous inequality constrained linear MPC problem is derived, in which the system dynamics are eliminated by using differential flatness. The problem is then converted into a polynomial optimization problem by parameterizing the flat output with piece-wise polynomials. 
		\item Path constraints are rigorously guaranteed by using sum of squares theory to transcribe the non-negative constraints of polynomials into the equivalent positive semi-definite constraints of matrices, and an equivalent SDP programming problem is formulated.
		\item The SDP problem is solved by using the customized primal-dual splitting-based iterations and accelerated by parallel computing.
	\end{itemize}

	It is worth noting that, although the derivations in this paper are carried out for linear MPC problems, it can also be extended to nonlinear MPC problems if the constraints remain linear and the objective remains quadratic after the differentially flat transformation.

	The paper is organized as follows. Differential flatness theory is stated and the form of flatness map for linear systems is described in Section \uppercase\expandafter{\romannumeral2}. The transformation of linear MPC problem with continuous-time path constraints into SDP problem is discussed in Section \uppercase\expandafter{\romannumeral3}. In Section \uppercase\expandafter{\romannumeral4}, we present the PDHG algorithm for SDP solving and explain that it can be accelerated by parallel computing. The simulation validation of our proposed MPC solver on quadruple-tank process is provided in Section \uppercase\expandafter{\romannumeral5}. Finally, concluding remarks are made in Section \uppercase\expandafter{\romannumeral6}

	
	
	
	
	%
	%
	%
	%
	%
	%


	\section{Preliminary: Differential Flatness of Linear System}\label{sec:flatness}
	
	Differential flatness is an important concept for a class of linear and nonlinear systems~\cite{763209}. A system is differentially flat if and only if there exists a flat output, such that all states and inputs subject to system dynamical constraints can be explicitly expressed as functions of the flat output (which is free of dynamical constraints) and a finite number of its derivatives.
	
	In this paper, we restrict our discussion to linear system. Consider an LTI system governed by the following ordinary differential equation:
	\begin{align}
		\dot{x}(t) = A x(t) + B u(t), \label{eq:sys}
	\end{align}
	where $x \in \mathbb{R}^n, u \in \mathbb{R}^m$. Without loss of generality, we assume that $(A,B)$ is controllable. Otherwise, we can always perform a Kalman decomposition and only consider the controllable part of the system.
	
	For such a system, Filess et al.\cite{763209} proved the following theorem:
	\begin{theorem}[Linear flatness~\cite{763209}]\label{th:linear_flat}
		A linear system is differentially flat, if and only if, it is controllable.
	\end{theorem}
	
	In general, the choice of the flat output may not be unique. In this paper, we adopt the procedure proposed by Yong et al.~\cite{7171938} to derive the flat output as well as the flatness map:
	\begin{theorem}
		If $(A,B)$ is controllable, then there exists a matrix $\mathcal{T}\in\Rb^{m\times n}$, such that the following $y$ is the flat output of the system:
		\begin{equation}
			y =\begin{bmatrix}
				y_1 & \cdots & y_m
			\end{bmatrix}^\top\triangleq \mathcal{T}x \in\Rb^m .	
		\end{equation}
		Moreover, there exists matrices $S\in\Rb^{n\times(n+m)}$, and $H\in\Rb^{m\times(n+m)}$, such that the state and the input of the system can be represented by the following flatness map:
		\begin{align*}
			x = S\bm{y}, u = H\bm{y},
		\end{align*}
		where $\bm{y}$ is the extended flat output vector consisting of $y_i$s and their derivatives, i.e.,\footnote{$\kappa_i,i\in\{1,\cdots,m\}$ is determined by $A,B$ and satisfies $\sum_{i=1}^{m} \kappa_i =n$.}
		{\small
			\begin{align}\label{eq:by}
				\bm{y}&\triangleq \begin{bmatrix} y_1 & \cdots & y_1^{(\kappa_1+1)} & \cdots & y_m & \cdots & y_m^{(\kappa_m+1)} \end{bmatrix}^\top \in\Rb^{n+m}.
			\end{align}
		}
	\end{theorem}
	\noindent
	The procedure to construct the $\mathcal{T},\,S,\,H$ matrices and extended flat output $\bm{y}$ (including the calculation of $\kappa_i$) is omitted due to space limit and the readers can refer to \cite{7171938} for more details.	
	
	
	\section{SDP Formulation of MPC}
	
\begin{figure}
		\centering
%
%

\begin{tikzpicture}

\node (problem1) at (0, 2)  [box, align=center] {Problem 1};
\node (problem3) at (0, 0) [box, align=center] {Problem 2};
\draw[<->,line width=2pt] (problem1) -- node[left, align=center,] {\small Differential \\  flatness} (problem3);

\node (problem4) at (5, 0) [box, align=center] {Problem 3};
\node (problem5) [box, align=center, above=1cm of problem4] {Problem 4};

\draw[->,line width=2pt] (problem3) -- node[left, align=center, right= -1.2cm of problem3] {\small Polynomial \\  \small parametrization} (problem4);
\draw[<->,line width=2pt] (problem4) -- node [left, align=right]{\small Markov-Luk\'{a}cs \\ \small Theorem} (problem5);

\end{tikzpicture}

		\caption{The relationships between optimization problems in this paper. Double tail arrow represents that the two problems are equivalent. Single tail arrow means that Problem \ref{pb:poly_opt} is obtained by parameterizing Problem \ref{pb:MPC_flat} using polynomials.}
		\label{fig:procedure}
	\end{figure}
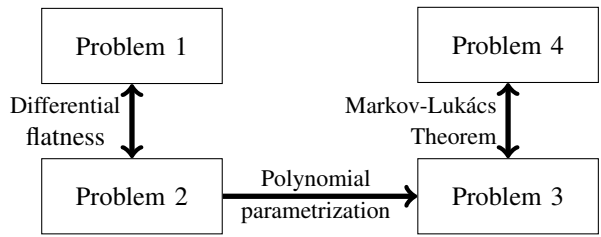

	This section is devoted to transcribing the MPC problem of a continuous-time path constrained linear system \eqref{eq:sys} to an SDP problem, the procedure of which is depicted in Fig~\ref{fig:procedure}. In the next subsection, we first remove the differential equality constraints in the MPC problem by differential flatness and then convert the problem into a polynomial optimization problem by parameterizing the flat output with piecewise polynomials. The polynomial optimization problem is then transformed into an equivalent SDP problem via Markov-Luk\'{a}cs Theorem~\cite{Roh_SOS,szego_orth_poly} and Sum-of-Squares (SOS) in Section~\ref{sec:sdpformulation}.
	
	\subsection{Polynomial Optimization Formulation of MPC}
	
	We consider an optimal control problem of the continuous-time linear system \eqref{eq:sys} under state and input constraints, which can be formulated in a receding horizon fashion as follows:  
	
	\begin{problem}[Continuous-time Linear MPC Problem]\label{pb:continuous}
		\begin{align*}
			\min_{x(t), u(t)} \quad & \int_0^T x(t)^\top Q x(t) + u(t)^\top R u(t)dt \\
			\text{s.t.} \quad & \dot{x}(t) = A x(t) + B u(t), \quad \forall t \in [0,T]\\
			&  \Xi x(t) + \Upsilon u(t) \leq b, \quad \forall t \in [0,T]\\
			& x(0)=x_0,
		\end{align*}
		where $T$ is the horizon length and $\Xi\in \mathbb R^{p\times n},\, \Upsilon\in \mathbb R^{p\times m}$ are matrices and $b\in \mathbb R^p$ is a vector of proper dimensions.
	\end{problem}
	
	Adopting the flatness map in Section~\ref{sec:flatness}, we can express the state $x(t)$ and control input $u(t)$ using the extended flat output $\bm{y}(t)$ and hence removing the differential equation constraint in Problem~\ref{pb:continuous}, which results in the following problem:  
	\begin{problem}[MPC using Flat Output]\label{pb:MPC_flat}
		\begin{align*}
			\min_{y(t)} \quad & \int_{0}^{T} \bm{y}(t)^\top (S^\top Q S+H^\top R H) \bm{y}(t) \ {\rm d}t \\
			\text{s.t.} \quad &  (\Xi S + \Upsilon H)\bm{y}(t) \leq b, \quad \forall t \in [0,T]\\
			&S\bm{y}(0)=x_0.
		\end{align*} 
	\end{problem}
	Notice that Problem~\ref{pb:continuous} and Problem~\ref{pb:MPC_flat} are equivalent, in the sense that we can use the definition of the flat output $y=\mathcal{T}x$ and the flatness map $x = S\bm{y}$, $u=H\bm{y}$ to map the solution of one problem to the other.
	
	Further notice that the path constraint $\Xi x(t) + \Upsilon u(t) \leq b$ (or $(\Xi S + \Upsilon H)\bm{y}(t) \leq b$), which consists of $p$ linear inequalities, requires that the state and the control input (or the flat output) to be inside a polytope at all time interval $[0,T]$. Aside from very special cases, Problem~\ref{pb:continuous} (or Problem~\ref{pb:MPC_flat}) cannot be solved in the infinite dimensional function space, due to the difficulty to determine when the path constraints are active~\cite{garcia1989model}.
	
	To facilitate optimization-based method to solve Problem~\ref{pb:MPC_flat}, we propose to parameterize the flat output $y(t)$ by piecewise polynomials, which effectively reduce the domain of the optimization problem from infinite dimensional function space to a finite dimensional space. To this end, first define the polynomial basis of degree $d$ as 
	\begin{align}\label{eq:defbeta}
		\gamma(t) = \begin{bmatrix}
			t^{d}& \cdots& t &1
		\end{bmatrix}^\top.
	\end{align}
	Suppose each entry of flat output $y$ is represented by $N$ segments of polynomials in the horizon $[0,T]$.
	Denote row vector $c_{l,i}\in\Rb^{(d+1)\times 1}$ as the coefficient of segment $l$ of flat output $y_i$, i.e.,\footnote{Since we need smoothness constraints on the conjecture points of segments, $c_{l,i},\cdots,c_{l+1,i}$ are not fully free and coupled by equality constraints.}
	\begin{align}\label{eq:by=cbeta}
		y_i(t)=\begin{cases}
			c_{1,i}^\top \gamma\left(\frac{tN}{T}\right), 0\leq t< \frac{T}{N}\\
			c_{2,i}^\top \gamma\left(\frac{tN}{T}-1\right), \frac{T}{N}\leq t< \frac{2T}{N} \\
			\hspace{60pt} \vdots \\
			c_{N,i}^\top \gamma\left(\frac{tN}{T}-(N-1)\right), \frac{(N-1)T}{N}\leq t< T
		\end{cases}.
	\end{align}
	Each segment of polynomial $c_{l,i}^\top \gamma\left(\cdot\right),l\in\{1,\cdots,N\}$ has been normalized such that the time variable is on interval $[0,1]$.
	
	By stacking the coefficients of the $l$-th segment $c_{l,i}$ vertically, we have the overall coefficient vector
	\begin{align}
		c_{l}\triangleq\begin{bmatrix}
			c_{l,1}\\ \vdots\\ c_{l,m}
		\end{bmatrix}\in\Rb^{m(d+1)\times 1}, \,
		\bc\triangleq\begin{bmatrix}
			c_1\\ \vdots\\ c_{N}
		\end{bmatrix}\in\Rb^{m(d+1)N\times 1}.
	\end{align}
	
	
	
	As a result, instead of optimizing $y$ in the infinite-dimensional function space, we can restrict ourselves to the following polynomial optimization problem:
	
	\begin{problem}[Polynomial Optimization]\label{pb:poly_opt}
		\begin{align*}
			\min_{\bm{c}}\ & J(\bm{c})= \bm{c}^\top P \bm{c}\\
			\text{s.t. }& (L_{j} c_{l}-g_{j})^\top \gamma(t) \geq 0, \forall t \in [0,1], \\
			&\qquad \qquad j\in\{1,\cdots,p\}, l\in\{1,\cdots,N\}\\
			&h_j^\top \bc=r_j\ ,j\in\{1,\cdots,2mN \}
		\end{align*}
	\end{problem}
	The calculation of parameters in Problem \ref{pb:poly_opt} is as follows.
	Define $\mathbf{e}_{j}$ as the canonical basis vector of length $d+1$, where 1 in on the $j$-th entry and $0$ on other entries.
	Define a matrix to represent the derivative of $d$ degree polynomial:
	\begin{align}
		D\triangleq \begin{bmatrix}
			0 &  &  &  &\\
			d & 0&  &  & \\
			& d-1 & 0 &  &\\
			&&\ddots&\ddots&\\
			& &  & 1 & 0 
		\end{bmatrix}\in\Rb^{(d+1)\times (d+1)}.
	\end{align} 
	Thus, for any coefficient $c\in\Rb^{d+1}$, we have polynomial derivative equation
	$\frac{{\rm d} \left(c^\top \gamma (t) \right)}{{\rm d} t}=(Dc)^\top \gamma(t).$
	Define $
	\bm{D}_{k}=\begin{bmatrix}
		I & D^\top & (D^2)^\top & \cdots & (D^k)^{\top}
	\end{bmatrix}^\top.
	$
	Then based on \eqref{eq:by}, one can verify that after polynomial parameterizing, the relationship between flat output $y$ and extend flat output $\bm{y}$ are $\bm{y}=\bm{\Pi} \cdot y$ with $\bm{\Pi}$ defined as
	\begin{align*}
		\bm{\Pi}\triangleq \begin{bmatrix}
			\bm{D}_{\kappa_1+1} &&&\\
			&	\bm{D}_{\kappa_2+1} &&&\\
			&	&\ddots &\\
			&	&& \bm{D}_{\kappa_m+1}
		\end{bmatrix}.
	\end{align*}
	Denote the $j$-th row of $\Xi,\Upsilon$ as $[\Xi]_j,[\Upsilon]_j$ respectively. Denote the $j$-th entry of vector $b$ as $b_j$.
	Then the parameters in Problem \ref{pb:poly_opt} is defined as:
	\begin{align*}
		L_{j}\triangleq &\left(([\Xi]_j S+ [\Upsilon]_j H)\otimes I_{d+1} \right) \times \bm{\Pi},\\
		g_j\triangleq & b_j \ \mathbf{e}_{d+1},
	\end{align*}
	where $\otimes$ is the Kronecker product and $\bm{0}_d$ is the all-zero vector of length $d$. $I_{d+1}$ is the identity matrix of size $(d+1)\times (d+1)$.
	
	The equality constraints are composed of segment smooth conditions and initial conditions, that is, for neighboring polynomial segments, the value of the polynomial and the value of its first-order derivative at conjecture points or at initial time are the same. 
	The equality constraint parameters are defined by 
	{\small
		\begin{align*}
			h_j=\begin{cases}
				\mathbf{e}^{N}_1\otimes \mathbf{e}_j \otimes \mathbf{e}_{d+1} , \text{ if } 1\leq j\leq m. \\
				\mathbf{e}^{N}_{l}\otimes \mathbf{e}_{j-m} \otimes \mathbf{e}_{d+1}-\mathbf{e}^{N}_{l+1}\otimes \mathbf{e}_{j-m} \otimes \mathbf{1}_{d+1} ,\\ \hspace{50pt}\text{ if } m+1\leq j\leq mN  .\\
				\mathbf{e}^{N}_1\otimes \mathbf{e}_j \otimes \mathbf{e}_{d} , \text{ if } mN+1\leq j\leq mN+m. \\
				\mathbf{e}^{N}_{l}\otimes \mathbf{e}_{j-mN-m} \otimes \mathbf{e}_{d+1}-\mathbf{e}^{N}_{l+1}\otimes \mathbf{e}_{j-Nm-m} \otimes \mathbf{1}_{d+1}, \\ \hspace{50pt}\text{ if } mN+1\leq j\leq mN+m .
			\end{cases}
		\end{align*}
	}\noindent
	where $\mathbf{e}^{N}_j$ is the canonical basis vector of size $N$, with 1 one $j$-th entry and 0 on other entries.
	\begin{align*}
		r_j=\begin{cases}
			[\Tc x_0]_j , \text{ if } 1\leq j\leq m .\\
			0, \text{ if } m+1\leq j\leq mN .\\
			[\Tc Ax_0]_{j-mN} , \text{ if } mN+1\leq j\leq mN+m .\\
			0 , \text{ if } mN+1\leq j\leq mN+m .
		\end{cases}
	\end{align*}
	
	Define matrix $P_{\rm int}\in \Rb^{(d+1)\times (d+1)}$ associating objective function integration on $[0,1]$ as
	\begin{align}
		[P_{\rm int}]_{u,v}=\frac{1}{2d+3-u-v}.
	\end{align}
	The parameter $P\in\Rb^{m{N}(d+1)\times mN(d+1)}$ in the objective of Problem \ref{pb:poly_opt} is calculated by
	\begin{align}
		P=\left((S\bm{\Pi})^\top Q S\bm{\Pi}+(H\bm{\Pi})^\top R H\bm{\Pi} \right)\otimes I_N \otimes  P_{\rm int}.
	\end{align}
	
	\begin{remark}
		Piecewise polynomials are chosen to represent the flat output for the following reasons:
		\begin{itemize}
			\item The set of polynomials are closed under derivative operation and is dense in the function space, as is shown by the Stone-Weierstrass theorem. Hence, we can approximate any continuous functions to arbitrary precision. In fact, one can also use polynomials to approximate the derivatives and high order derivatives of a smooth enough function\cite{4908942}.
			\item The continuous-time path constraints are transformed into non-negativity of a univariate polynomial inside an interval, which can be transformed {\bf exactly} into Positive Semi-Definite (PSD) cone constraint using Markov-Luk\'{a}cs theorem and SOS \cite{Roh_SOS}. The detailed discussion is reported in the subsequent subsection.
		\end{itemize}
	\end{remark}
	
	\subsection{SDP Formulation via SOS}
	\label{sec:sdpformulation}
	
	This subsection is devoted to the \emph{exact} SDP formulation of the polynomial optimization Problem~\ref{pb:poly_opt}. To this end, the following theorem is needed:
	
	\begin{theorem}[Markov-Luk\'{a}cs theorem \cite{Roh_SOS}]
		Let $a<b$. Then, a polynomial $p(t)$ is non-negative for $t\in [a, b]$, if and only if it can be written as
		\begin{align*}
			p(t)=\begin{cases}
				f(t)+(t-a)(b-t) g(t), &  \text { if } \deg(p) \text { is even } \\
				(t-a) f(t)+(b-t) g(t), &  \text { if } \deg(p) \text { is odd }
			\end{cases},
		\end{align*}
		where $f(t), g(t)$ are SOS polynomials, with degree $\deg(f) \leq \deg(p)$, $\deg(g) \leq \deg(p)-2$ when $\deg(p)$ is even, or $\deg(f) \leq \deg(p)-1$, $\deg(g) \leq \deg(p)-1$ when $\deg(p)$ is odd.
	\end{theorem}
	
	For simplicity, we shall only consider the case where the flat output $y$ is an odd degree polynomial, i.e., $d$ is an odd number. The case where $d$ is even can be treated similarly. Let us denote $\delta \triangleq \frac{d-1}{2}$. Notice that a degree $d-1$ SOS polynomial $f$ can be represented as
	\[
	f(t)= \tilde \gamma(t)^\top X\tilde \gamma(t),
	\]
	with $\tilde \gamma(t)\triangleq\begin{bmatrix}
		t^{\delta}& \cdots& t &1
	\end{bmatrix}^\top$ and positive semi-definite matrix $X\in\Rb^{(\delta+1)\times (\delta+1)}$.
	
	As a result, each inequality constraint $(L_{j}c_l-g_{j})^\top \gamma(t)\geq 0 $ in Problem \ref{pb:poly_opt} can be equivalently represented as
	\begin{align}
		&	(L_{j}c_l-g_{j})^\top\gamma(t) =\notag \\
		&\qquad t \tilde \gamma(t)^\top X^f_{j,l}\tilde \gamma(t)+\left(1-t\right)\tilde \gamma(t)^\top X^g_{j,l}\tilde \gamma(t),\label{eq:odd_eq}
	\end{align}
	with positive semi-define matrices $X^f_j, X^g_j\in\Rb^{(\delta+1)\times (\delta+1)}$. By comparing the coefficients of the polynomials on the LHS and RHS of \eqref{eq:odd_eq}, we know that \eqref{eq:odd_eq} is equivalent to:
	\begin{align}
		L_{j}c_l-g_{j}=\Mc( X^f_{j,l}, X^g_{j,l})
	\end{align}
	where 
	\begin{align}\label{eq:defMc}
		\Mc(X^f_{j,l}, X^g_{j,l})=\begin{bmatrix}
			\tr(F_0 X^f_{j,l})+\tr(G_0 X^g_{j,l})\\
			\tr(F_1 X^f_{j,l})+\tr(G_1 X^g_{j,l})\\
			\vdots\\
			\tr(F_{d} X^f_{j,l})+\tr(G_{d} X^g_{j,l})
		\end{bmatrix},
	\end{align}
	and $\{F_i,G_i|i=0,1,\cdots,d\}$ is a sequence of constant matrices defined as
	\begin{align*}
		[F_i]_{u,v}=&\begin{cases}
			1, &\text{ if } u+v=i+2 \\
			0, &\text{ otherwise }
		\end{cases},\\
		[G_i]_{u,v}=&\begin{cases}
			-1, &\text{ if } u+v=i+2 \\
			1, &\text{ if } u+v=i+1 \\
			0, &\text{ otherwise }
		\end{cases},
	\end{align*}
	where $[\cdot]_{u,v}$ represents the entry at row $u$, column $v$ in a matrix.
	Notice that linear function $\Mc(\cdot,\cdot)$ is independent of index $j,l$.
	

	Now we handle the second order objective function $J(c)=\bm{c}^\top P \bm{c}$ by linear matrix inequality techniques. Notice that $P$ is a positive semi-definite matrix, define $\tP\in\Rb^{\rank(P)\times \size(P)}$, such that $\tP^\top\tP=P$. 
	Notice that the following three optimization problems are equivalent where $s$ is a scalar:
	\begin{align*}
		\min_{\bc\in\mathcal{C}}&\quad \bc^\top P \bc \Leftrightarrow
		\min_{\bc\in\mathcal{C},s}\ s, \text{ s.t. } \|\tP\bc \|_2\leq s\\
		\Leftrightarrow&\min_{\bc\in\mathcal{C},s\geq 0}s ,\text{ s.t. } \begin{bmatrix}\tP\bc  \\ s \end{bmatrix}\in\text{second order cone}.
	\end{align*}
	We arrive at the following SDP problem which is equivalent to Problem \ref{pb:poly_opt} and computationally tractable.
	\begin{problem}[SDP Problem]\label{pb:conic}\hspace{40pt}	
		
		\noindent
		\textbf{Original form}
		\begin{align}
			&\min_{\bm{c},s,\{X^f_{j,l},X^g_{j,l}\}  }  \quad   s \notag \\
			\text{s.t. }& L_{j} c_l-g_j=\Mc(X^f_{j,l},X^g_{j,l}),  X^f_{j,l},\,X^g_{j,l}\in \mathbb S_+,\notag \\
			& \ \qquad \qquad j\in\{1,\cdots,p \},  l\in\{1,\cdots,N \}\label{eq:ineq_eq_cons}\\
			&h_j^\top \bc=r_j ,j\in\{1,\cdots,2mN \}\label{eq:eq_eq_cons}\\
			&\begin{bmatrix}\tP \bc\\ s \end{bmatrix} \in\SOC \notag
		\end{align}
		where $\SOC$ denotes the {second order cone} and $\Sb_{+}$ is the positive semi-definite cone.
	\end{problem}
	For notation conciseness, we define:
	\begin{align}
		\bX\triangleq \begin{bmatrix}
			X^f_{1,1}&&&&\\
			&X^g_{1,1}&&&\\
			&&\ddots&&\\
			&&&X^f_{p,N}&\\
			&&&&X^g_{p,N}\\
		\end{bmatrix}.
	\end{align} 
	
	Moreover, define function
	\begin{align}\label{eq:defbM}
		\bM(\bX)=\begin{bmatrix}
			\Mc(X_{1,1}^f,X_{1,1}^g)\\
			\vdots \\
			\Mc(X_{p,N}^f,X_{p,N}^g)
		\end{bmatrix}
		=\begin{bmatrix}
			\begin{bmatrix}
				\tr(M_{1,1}^0\bX)\\
				\vdots\\
				\tr(M_{1,1}^{d}\bX)
			\end{bmatrix}\\
			\vdots \\
			\begin{bmatrix}
				\tr(M_{p,N}^0\bX)\\
				\vdots\\
				\tr(M_{p,N}^{d}\bX)
			\end{bmatrix}
		\end{bmatrix},
	\end{align}
	where $M_{j,l}^i$ is the corresponding matrix composed of $F_{i},G_i$ at compatible position that generates $\Mc(X_{j,l}^f,X_{j,l}^g)$ in \eqref{eq:defMc}. 
	Define 
	\begin{align}
		\bL=& I_{N} \otimes \begin{bmatrix}
			L_1^\top & L_2^\top & \cdots & L_p^\top
		\end{bmatrix}^\top ,\\
		\bg=&I_{N} \otimes \begin{bmatrix}
			g_1^\top & g_2^\top & \cdots & g_p^\top
		\end{bmatrix}^\top ,\\
		\bh=& \begin{bmatrix}
			h_1^\top & h_2^\top & \cdots & h_{2mN}^\top
		\end{bmatrix}^\top, \\
		\br=& \begin{bmatrix}
			r_1^\top & r_2^\top & \cdots & r_{2mN}^\top
		\end{bmatrix}^\top.
	\end{align}
	We can rewrite Problem~\ref{pb:conic} as
	
	\textbf{\textit{Compact form}}
	\begin{align}
		&\min_{s,\bc,\bX}\quad s\notag \\
		\text{s.t. }
		&\bL \bc-\bM(\bX)=\bg \label{eq:ineq_comp}\\
		&\bh \bc=\br \label{eq:eq_comp} \\
		&\bX\in\Sb_+,s\geq 0 \notag \\
		&\begin{bmatrix}\tP \bc \\ s \end{bmatrix} \in\SOC \notag
	\end{align}
	
	Since there are $p$ inequality constraints in the original Problem~\ref{pb:continuous}, $\bm{X}$ is a block diagonal matrix with $2pN$ positive semi-definite matrices of size $\delta+1$. As a result, in the following section, we introduce a customized algorithm that solves Problem \ref{pb:conic} by primal-dual hybrid gradient methods which can handle $\bX$ in a parallel fashion. 
	However, before continuing on, we would like to give a comparison between the conventional quadratic programming-based linear MPC and our approach.
	
	\subsection{Discussions}
	
	A conventional way to solve the continuous-time optimal control problem is to discretize it into the following discrete-time linear MPC problem\cite{garcia1989model}:
	
	\begin{problem}[Discrete-time Linear MPC Problem]\label{pb:discrete}
		\begin{align*}
			\min_{ \{x[k],u[k]\}_{k=1}^{T_d} }&\quad  \sum_{k=1}^{T_d} x[k]^\top Q x[k] + u[k]^\top R u[k] \\
			\text{s.t.} \quad & x_{k+1} = A_d x_{k} + B_d u_k, \quad k = 1, \cdots, {T_d} \\
			& \Xi x[k] + \Upsilon u[k]  \leq b, \quad k = 1, \cdots, {T_d}
		\end{align*}
	\end{problem}
	\noindent  where $A_d,B_d$ are the discretized system matrix assuming zero-order hold for the control input is used and $T_d\in\Zb_{+}$ is the discrete horizon length.
	
	One of the main differences between Problem~\ref{pb:discrete} and Problem~\ref{pb:poly_opt} is that the control input is parameterized as step functions in Problem~\ref{pb:discrete} (assuming zero-order hold is used), while for our case, the flat output (and hence the control input as it is a linear function of the flat input and its derivatives) is parameterized as polynomials. 
	
	Another difference is that Problem~\ref{pb:discrete} is a Quadratic Programming (QP) problem and hence can be solved more efficiently than SDP. However, this is due to the fact that in Problem~\ref{pb:discrete}, constraints are only required to hold at discrete sampling time instants and therefore they may be violated in between sampling times. 
	
	On the other hand, the reason for our SDP formulation is that we want to have an exact representation of the continuous time path constraints. If we only require the constraint to hold at discrete time instant, since the value of a polynomial at a time instant is a linear function of its coefficients, we can express such constraints as linear inequalities on the coefficients of the polynomial, which effectively relaxed the polynomial optimization Problem~\ref{pb:poly_opt} into a QP problem that only guarantees constraint satisfaction at a discrete time instant. As an alternative, one could also leverage the following theorem to generate a QP problem, which has a smaller feasible set than that of Problem~\ref{pb:poly_opt}, but is guaranteed to satisfy the path constraints at every time instant. 
	
	\begin{theorem}[\cite{powers2000polynomials}]
		Let ${\rm Pd}([a,b])$ denote the set of polynomials $p(t) > 0, \forall t \in [a,b]$. Define 
		\begin{align*}
			\mathcal{P}_q:=\left\{\left.\sum_{i+j \leq q} c_{ij}(b-x)^i(x-a)^j \right| c_{ij} \geq 0\right\}.
		\end{align*}
		If polynomial $p \in {\rm Pd}([a,b])$, then $p \in \mathcal{P}_q$ for sufficiently large integer $q$.
	\end{theorem}
	
	\section{Accelerated SDP Solving with Parallel Computing}
	
	\subsection{Primal dual hybrid gradient for SDP solving}

	
	In this subsection, we present the primal-dual hybrid gradient algorithm that solves Problem \ref{pb:conic}. 
	Encode the constraints into the objective function as
	\begin{align}
		\min_{s,\bc,\bX} s+
		\Ib_{\Sb_{+}, \SOC}(\bX,\tbc,s)+\Ib_{=}(\bX,\bc,\tbc), 
	\end{align}
	where $\tbc=\tP \bc$ is the slack variable, and the indicator functions are defined as 
	\begin{align}
		&	\Ib_{\Sb_{+}, \SOC}(\bX,\tbc,s)=\begin{cases}
			0,\text{ if }\bX\in\Sb_+ \text{ and}\begin{bmatrix}\tbc \\ s \end{bmatrix} \in\SOC\\
			+\infty, \text{otherwise}
		\end{cases}.\\
		&\Ib_{=}(\bX,\bc,\tbc)=\begin{cases}
			0,\text{ if }\left\{
			\makecell{
				\bL \bc-\bM(\bX)=\bg\\
				\bh \bc=\br\\
				\tP \bc-\tbc=0
			} \right. \\
			+\infty, \text{otherwise}
		\end{cases}.
	\end{align}
	
	Using primal-dual operator splitting \cite{exp_lowrank}, the iterations can be derived as the following where $\alpha$ is the primal step-size and $\beta$ is the dual step-size. $\Dc^*_{\bX}(\cdot)$ is the conjugate operator of the linear mapping $\bM(\bX)$, and $\Dc^*_{\bc}(\cdot)$ is the conjugate operator of the linear mapping $\begin{bmatrix}
		\bL^\top& \bh^\top  & \tilde{P}^\top \end{bmatrix}^\top\bc$, the definition of which are given in \eqref{eq:D*X} and \eqref{eq:D*c} respectively. 
	
	\textbf{1. Primal step:}
	\begin{align}
		\bX^{k+1} &\leftarrow \proj_{\Sb_+}(\bX^{k}-\alpha \Dc_{\bX}^*(\lambda_1^k))	\label{eq:primal_X}\\
		\bc^{k+1}&\leftarrow \bc^{k}-\alpha \Dc_{\bc}^*(\lambda_1^k,\lambda_2^k,\lambda_3^k)\label{eq:primal_c}\\
		\begin{bmatrix}
			\tbc^{k+1}\\
			s^{k+1}
		\end{bmatrix}&\leftarrow \proj_{\SOC}\begin{bmatrix}
			\tbc^{k}-\alpha (-\lambda_3^k)\\
			(s^{k}-\alpha)^{+} \label{eq:primal_soc}
		\end{bmatrix}
	\end{align}
	where $(s^{k}-\alpha)^+=\max(s^{k}-\alpha,0)$.
	
	\textbf{2. Calculating difference:}
	\begin{align}
		\Delta \bX^{k+1}&\leftarrow 2\bX^{k+1}-\bX^{k}\label{eq:diff_X}\\
		\Delta \bc^{k+1}&\leftarrow 2\bc^{k+1}-\bc^{k}\\
		\Delta \tbc^{k+1}&\leftarrow 2\tbc^{k+1}-\tbc^{k}
	\end{align}
	
	\textbf{3. Dual step:}
	\begin{align}
		\lambda_1^{k+1}&\leftarrow 	\lambda_1^k+\beta \bL \Delta\bc^{k+1} -\beta\bM(\Delta\bX^{k+1})-\beta\bg  \label{eq:dualy1}\\
		\lambda_2^{k+1}&\leftarrow 	\lambda_2^k+\beta \bh \Delta\bc^{k+1}-\beta\br\\
		\lambda_3^{k+1}&\leftarrow 	\lambda_3^k+\beta \tP \Delta\bc^{k+1}-\beta\Delta \tbc^{k+1}\label{eq:dualy3}
	\end{align}
	Here $\lambda_1$ is the dual variable corresponding to the equality constraint \eqref{eq:ineq_comp} and equivalently \eqref{eq:ineq_eq_cons}. $\lambda_2$ is the dual variable corresponding to the equality constraint \eqref{eq:eq_comp} and equivalently \eqref{eq:eq_eq_cons}. $\lambda_3$ is the dual variable corresponding to the equality constraint $\tP \bc-\tbc=0$.
	Define $\lambda=\begin{bmatrix}
			\lambda_1^\top & \lambda_2^\top & \lambda_3^\top 
		\end{bmatrix}^\top$.
		
		We denote the entry of $\lambda_1$ corresponding to segment $l$, inequality index $j$, order $i$ as $\lambda_1[i,j,l],i\in\{0,\cdots,d\}, j\in\{1,\cdots,p\},l\in\{1,2,\cdots,N\}$. Recall the definition of $M_{j,l}^i$ in \eqref{eq:defbM}. The conjugate operators $\Dc^*$ are defined as\footnote{ We use $[A]_j$ to denote the $j$-th row of matrix $A$, or in case of $A$ is column vector, $j$-th entry of the vector $A$.}
		
		\begin{align}
			\Dc^*_\bX(\lambda_1)\triangleq&-\sum_{l=1}^{N}\sum_{j=1}^{p} \sum_{i=1}^{d+1} \lambda_1[i,j,l] M^i_{j,l} \label{eq:D*X}\\
			\Dc^*_\bc(\lambda_1,\lambda_2,\lambda_3)\triangleq& \sum_{l=1}^{N}\sum_{j=1}^{p} \sum_{i=1}^{d+1} \lambda_1[i,j,l]\left(\mathbf{e}^N_l\otimes [L_{j}]_i^\top  \right)\notag\\
			&+\sum_{j=1}^{2mN} [\lambda_2]_j h_j+ \tP^\top \lambda_3 \label{eq:D*c}
		\end{align}
		where $\mathbf{e}^N_l$ is the canonical basis vector of size $N$, with 1 on $l$-th entry and $0$ on other entries.
		
		The projection to the semi-definite cone is 
		\begin{align}
			\proj_{\Sb_+}(X)=\sum_{i=1}^{\size(X)} \max \left\{0, \nu_i\right\} \mu_i \mu_i^\top
		\end{align}
		where $\nu_i,\mu_i$ are the eigenvalue and the corresponding eigenvector of $X$. The projection to the second order cone is 
		\begin{align}
			\proj_{\SOC}\begin{bmatrix}
				\bc\\
				s
			\end{bmatrix}=
			\begin{cases}
				\frac{s+\|\bc\|_2}{2\|\bc\|_2}\begin{bmatrix}
					\bc\\
					\|\bc\|_2
				\end{bmatrix}& \text{if } \|\bc\|_2>s.\\
				\begin{bmatrix}
					\bc\\
					s
				\end{bmatrix}& \text{if } \|\bc\|_2\leq s.
			\end{cases}.
		\end{align}
		The convergence of algorithm \eqref{eq:primal_X}-\eqref{eq:dualy3} is provided in the following.
		\begin{theorem}[\cite{ryu_yin_2022}]
			Assume the solution to KKT conditions of Problem \ref{pb:conic} exists (denoted by $\bc^\star,\bX^*,s^*,\lambda^*$), and strong duality holds. If the linear projection defined by 
			\begin{align*}
				\Lc(\bc,\bX)=\begin{bmatrix}
					\bL\bc-\bM(\bX) \\ \bh\bc 
				\end{bmatrix}
			\end{align*} 
			and step sizes $\alpha,\beta$ satisfy $0<\alpha\beta < 1/\left\| \Lc(\bc,\bX) \right\|_2$, then the primal dual hybrid gradient descent algorithm \eqref{eq:primal_X}-\eqref{eq:dualy3} converges to the solution to KKT conditions, i.e., $\bc^k\rightarrow \bc^\star,\bX^k\rightarrow \bX^\star,s^k\rightarrow s^\star,\lambda^k\rightarrow \lambda^\star$.
		\end{theorem}
		
		\subsection{GPU parallel computing}\label{subsec:GPU}
		
		It is worth noticing that for our proposed iterations, a significant proportion of the time will be spent on the projection $\proj_{\Sb_+}(\cdot)$. However, since $\bX$ is a block diagonal matrix with $2pN$ matrices of size $\delta +1$ on its diagonal, the projection of $\bX$ can be parallelized by projecting each small matrices onto the PSD cone. Furthermore, the calculation of $\Dc^*_{\bX},\Dc^*_{\bc}$ in \eqref{eq:D*X} and \eqref{eq:D*c}, and the difference calculation in \eqref{eq:diff_X} are essentially tensor operations and hence can be accelerated by parallel computation.
		
		To speed up the computation of the proposed PDHG solver, we implement it in a parallelized manner on the GPU. Specifically, the projection step~\eqref{eq:primal_X} is wrapped as a kernel to be computed in parallel on the GPU. Additionally, we implement the calculation of $\mathcal{D}^*_{\boldsymbol{X}}$ and $\mathcal{D}^*_{\boldsymbol{c}}$ in~\eqref{eq:D*X} and~\eqref{eq:D*c}, as well as the update steps from~\eqref{eq:primal_c} to~\eqref{eq:dualy3} as tensor operations, which can also be accelerated by GPU parallelization.
		
		We test the implemented our proposed solver on a desktop computer equipped with an AMD Ryzen Threadripper 3970X 32-Core Processor and an NVIDIA GeForce RTX 3080 GPU. We report the computation time for a single iteration of the parallelized solver running on GPU in Table~\ref{tab:pdhg_time}, and compare it to that of a serialized version running on the CPU. Our results show that the iteration time of the accelerated solver is significantly shorter than that of the CPU version. Moreover, the computation time of the accelerated solver increases slowly as the problem size (the value $d$ and $N$) grows. Even for the largest problem instances considered, the iteration time remains within a few milliseconds, demonstrating the effectiveness of the GPU acceleration and the efficiency of the implementation.
		\begin{table}[!htbp]
			\centering
			\begin{tabular}{lllll}
				\toprule
				\multirow{2}{*}{$N$} & \multicolumn{3}{c}{$d$} \\
				\cmidrule(lr){2-4}
				& \multicolumn{1}{c}{3}  & \multicolumn{1}{c}{5}  & \multicolumn{1}{c}{7} \\
				\midrule
				2 & 0.318 (1.874) & 0.345 (2.773) & 0.390 (3.778) \\
				6 & 0.408 (5.159) & 0.443 (8.012) & 0.435 (11.240) \\
				10 & 0.431 (8.811) & 0.506 (13.758) & 0.547 (19.254) \\
				20 & 0.484 (19.449) & 0.649 (30.721) & 0.821 (42.886) \\
				30 & 0.561 (32.330) & 0.767 (50.046) & 1.187 (63.170) \\
				40 & 0.776 (46.167) & 1.394 (65.989) & 2.593 (93.495) \\
				\bottomrule
			\end{tabular}
			\caption{Time consumption (in milliseconds) of a single iteration of the PDHG solver (calculating updates \eqref{eq:primal_X}-\eqref{eq:dualy3}) for different problem sizes. The runtime on CPU is in parentheses. 
			}
			\label{tab:pdhg_time}
		\end{table}
		
		
		%
		%
		
		\subsection{Warm Start and Termination Rule}
		
		Denote $\Delta t$ as the control apply time length of each solution. To facilitate the simple warm start strategy, the horizon $T$ satisfies $T=N\cdot \Delta t$, i.e., the first segment of control input $u(t)$ is applied before receding to a new horizon. We evaluate our algorithm's performance in two different strategies: cold start, and warm start.
		The cold start strategy initializes the optimization variable $\bc(t+\Delta t),\bX(t+\Delta t),\lambda(t+\Delta t)$ as random vectors/matrices with each entry uniformly distributed on $[-0.5,0.5]$. The warm start strategy initializes the optimization variable in a shifting manner, i.e., $\forall \ 2\leq l\leq N,1\leq j\leq p$:
		\begin{align*}
			\left(X^f_{j,l}(t+\Delta t)\right)^{0}&:=\left(X^f_{j,l-1}(t)\right)^{\tau(t)}\\
			\left(X^g_{j,l}(t+\Delta t)\right)^{0}&:=\left(X^g_{j,l-1}(t)\right)^{\tau(t)}\\
			\left(c_l(t+\Delta t)\right)^{0}&:=\left( c_{l-1}(t) \right)^{\tau(t)}
		\end{align*}
		where the super script $\tau(t)$ is the number of iterations applied to solve the SDP problem at time $t$.
		For the first segment, $\left(X^f_{j,1}(t+\Delta t)\right)^{0},\left(X^g_{j,1}(t+\Delta t)\right)^{0},\left(c_1(t+\Delta t)\right)^{0}$ are initialized randomly.
		The dual variable $\lambda$ is also initialized in a shifting manner according to its correspondence with $\bX,\bc$ in \eqref{eq:dualy1}-\eqref{eq:dualy3}.
		
		The algorithm termination is determined by the residue $\epsilon^k=\epsilon_{\text {primal }}^k+\epsilon_{\mathrm{dual}}^k $ where the primal and dual residue are defined as:
		\begin{align*}
			\epsilon_{\mathrm{primal}}^k=&\left\|\frac{1}{\alpha}\left(\bX^k-\bX^{k-1}\right)-\mathcal{D}_{\bX}^*\left(\lambda^k-\lambda^{k-1}\right)\right\|_F \\
			& +\left\|\frac{1}{\alpha}\left(\bc^k-\bc^{k-1}\right)-\mathcal{D}_{\bc}^*\left(\lambda^k-\lambda^{k-1}\right)\right\|_2 
		\end{align*}
		\begin{align*}
			\epsilon_{\mathrm{dual}}^k=&\left\|\frac{1}{\beta}\left(\lambda^k-\lambda^{k-1}\right)\right.\\
			& \left.-\begin{bmatrix}
				\bL(\bc^k-\bc^{k-1})-\bM(\bX^k-\bX^{k-1})-\bg\\
				\bh(\bc^k-\bc^{k-1})-\br\\
				\tilde{P}(\bc^k-\bc^{k-1}) - (\tbc^k-\tbc^{k-1})
			\end{bmatrix}\right\|_2.
		\end{align*}
		$\|\cdot\|_F$ is the Frobenius norm.
		Our proposed algorithm terminates when $\epsilon^k$ corresponding to segment $1$ is below $2\times 10^{-2}$, which is accurate enough for control performance. 
		The control and computational speed performance is demonstrated in the next section.
		
		\section{Simulation}
\begin{figure}[!htbp]
	\centering 
	\includegraphics[width=0.8\linewidth]{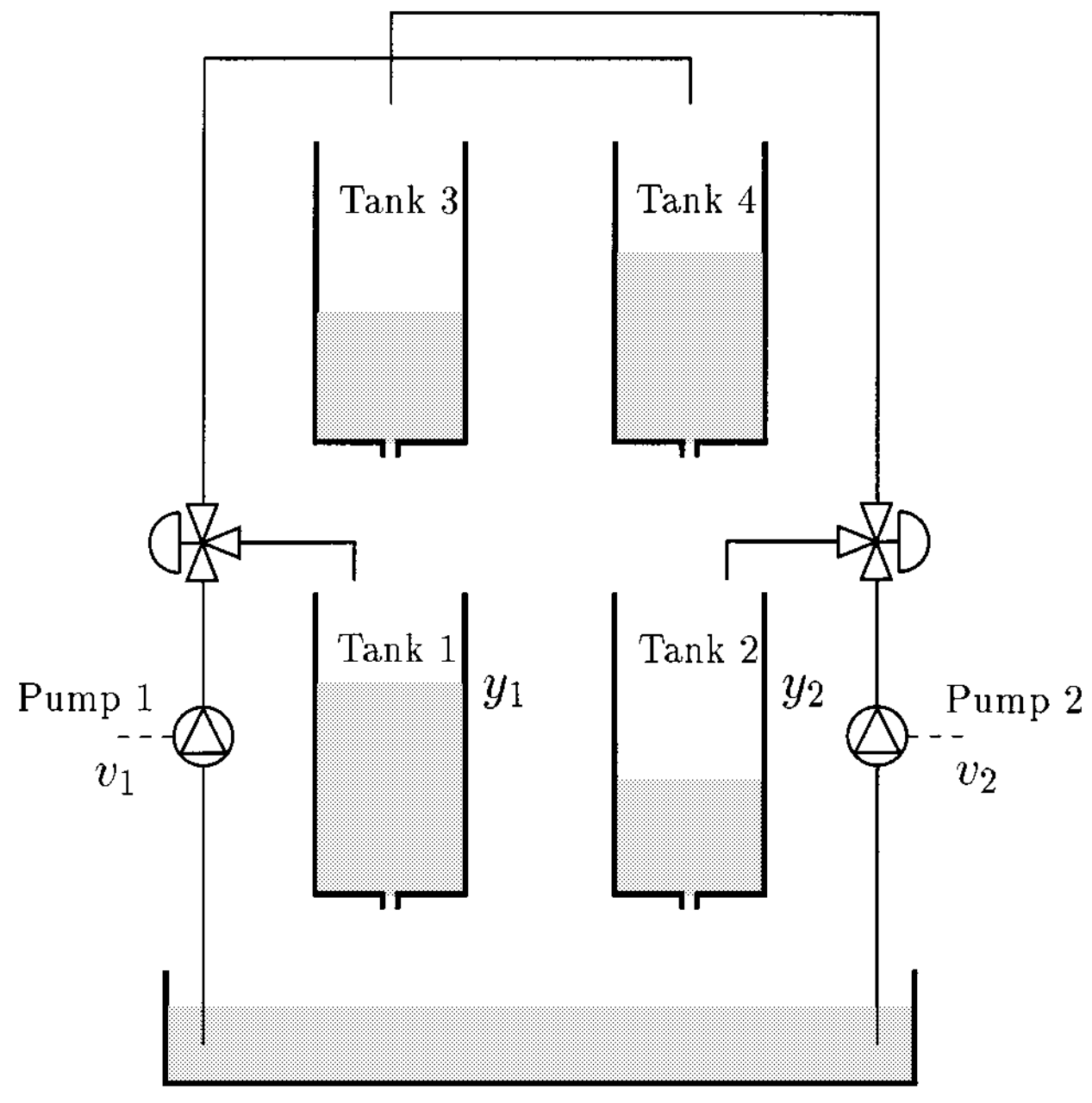}
	\caption{The schematic diagram of the quadruple-tank process.}
	\label{fig:quaduple-tank}
\end{figure}

		The performance of the proposed MPC solver is validated on the quadruple-tank process~\cite{johansson2000quadruple}, whose schematic diagram is visualized in Fig.~\ref{fig:quaduple-tank}. The system has $4$ states, which represent the liquid levels (in centimeter) of each tank. There are two control inputs in the system, namely the voltage (in volt) of Pump 1 and Pump 2 in Fig.~\ref{fig:quaduple-tank}. The simulation employs the same linearized system equations and system parameters as~\cite{johansson2000quadruple}, which are omitted due to space limitations.

The initial conditions of the tanks are 
\[x_0=\begin{bmatrix}
	10 & 19 & 19 & 1
\end{bmatrix}^\top,\]
and the control objective of the MPC is to track a reference trajectory $r(t)$ of the liquid levels. For simplicity, we set the constant reference trajectory as
\[r(t)=\begin{bmatrix}
	19.9 & 19.9 & 2.4 & 2.4 
\end{bmatrix}^\top.\]
In addition to tracking the reference signal, the MPC must ensure that the liquid levels in all tanks remain between 0 to 20 cm and that the control inputs stay in the voltage limit between 0 to 8 V during the control process. The objective weighting matrices in MPC Problem \ref{pb:continuous} are defined as $Q=I,R=0.1I$ with appropriate dimensions. Our code is available on  \href{https://github.com/zs-li/MPC_PDHG}{https://github.com/zs-li/MPC\_PDHG}.

\subsection{Control Performance}

\begin{figure}[!htbp]
	\vspace{-7pt}
	\centering
	\includegraphics[width=\linewidth]{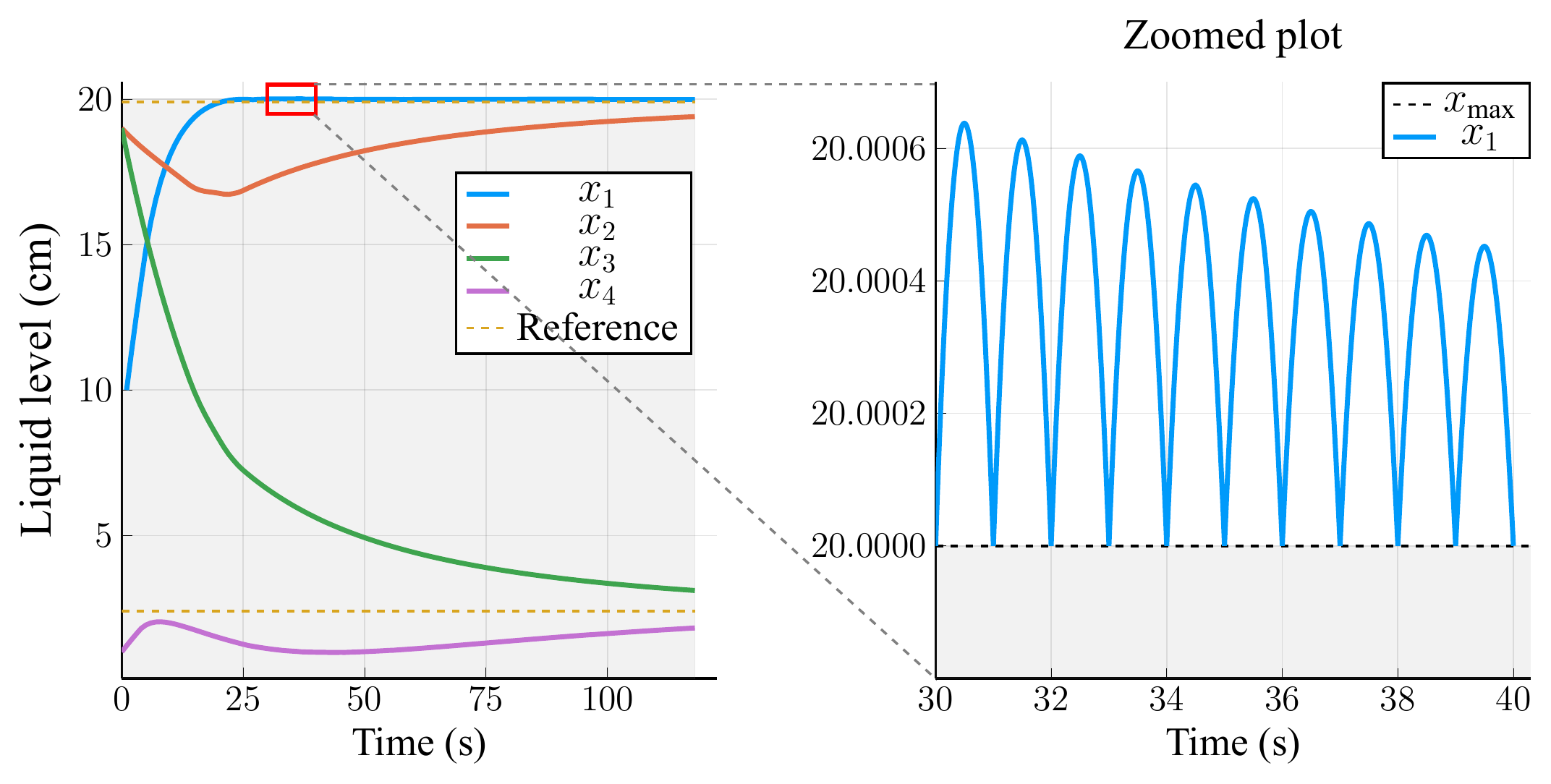}
	\caption{The states of discrete time linear MPC using the QP solver. The gray area in the figures denotes the feasible region of states. The states between sampling times violate the constraints.}\label{fig:QP}
	\vspace{-7pt}\end{figure}

\begin{figure}
	\vspace{-7pt}
	\centering
	\begin{subfigure}{0.45\linewidth}
		\includegraphics[width=\linewidth]{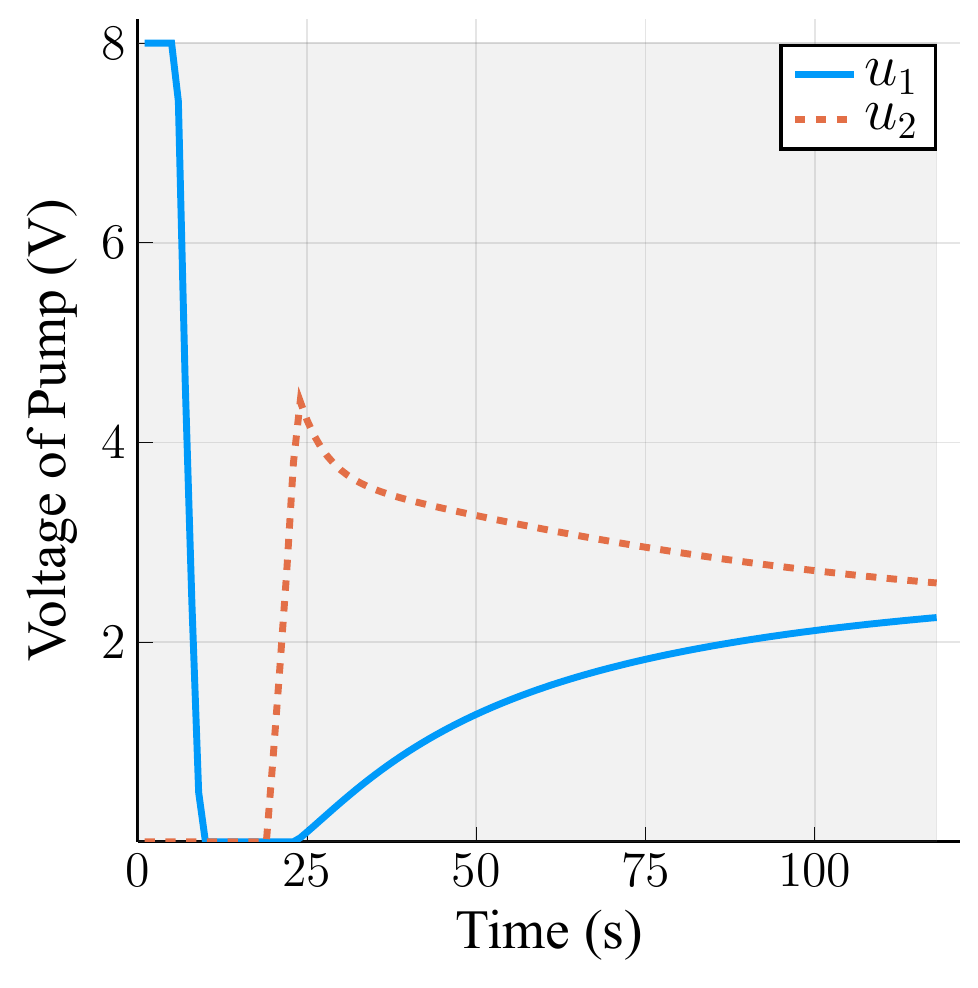}
		\caption{Control input using the QP solver.}
	\end{subfigure}
	\begin{subfigure}{0.45\linewidth}
		\includegraphics[width=\linewidth]{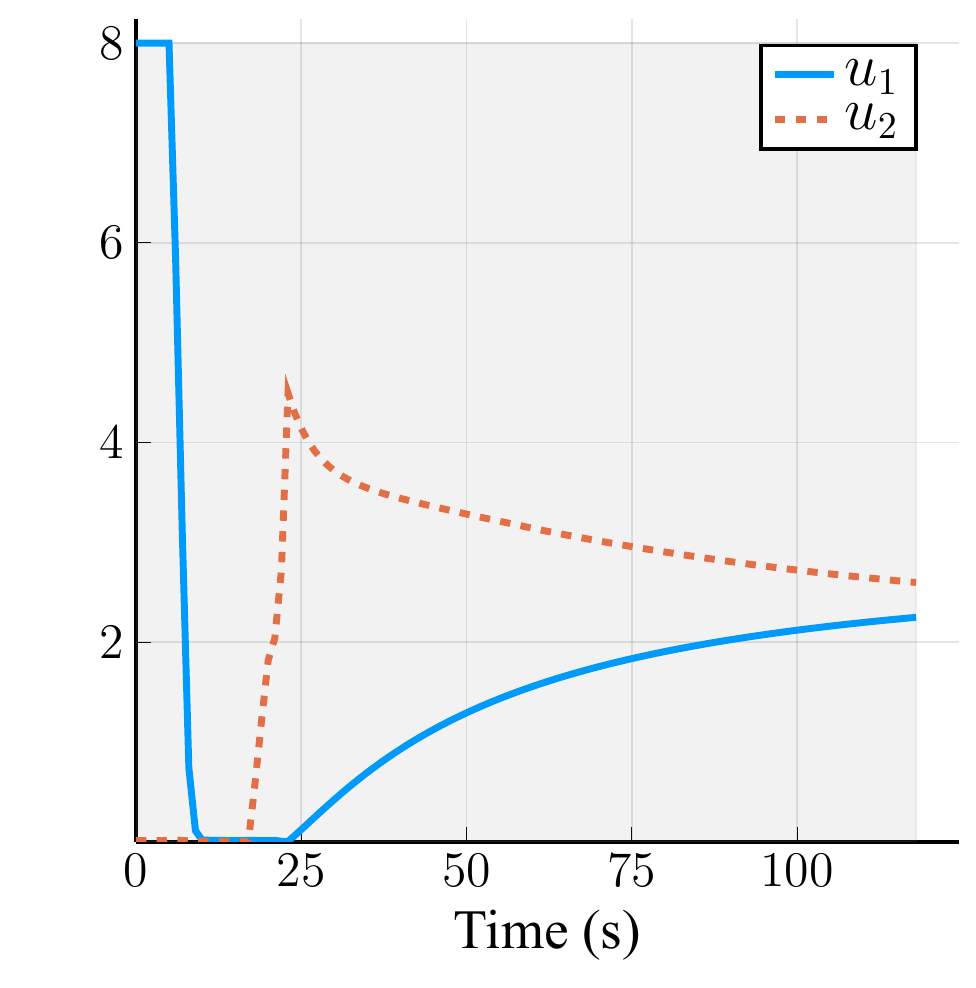}
		\caption{Control input using the proposed solver.}
	\end{subfigure}
	\caption{Comparison on the input trajectory using the QP and the proposed SDP strategy respectively. The gray area in the figures denotes the feasible regions of control input.}
	\label{fig:state_input}
	\vspace{-7pt}\end{figure}

For comparison, we employ the Quadratic Programming (QP) formulation (Problem~\ref{pb:discrete}), where the MPC problem is discretized with a sampling interval of $T_s=1{\rm s}$ and horizon length $T_d=20$. On the other hand, for the proposed method, we set the degree of polynomial $d=3$ and the segments of polynomials $N=20$, horizon length $T=20$. For each iteration, the resulting control input applies  to the system for $\Delta t=1$ second.
Thus, the two methods are comparable in terms of horizon length and update frequency. We simulate the control process for $120$ seconds and visualize the resulting system states and control inputs in Fig.~\ref{fig:QP}-Fig.~\ref{fig:state_input}.

\begin{figure}
	\includegraphics[width=\linewidth]{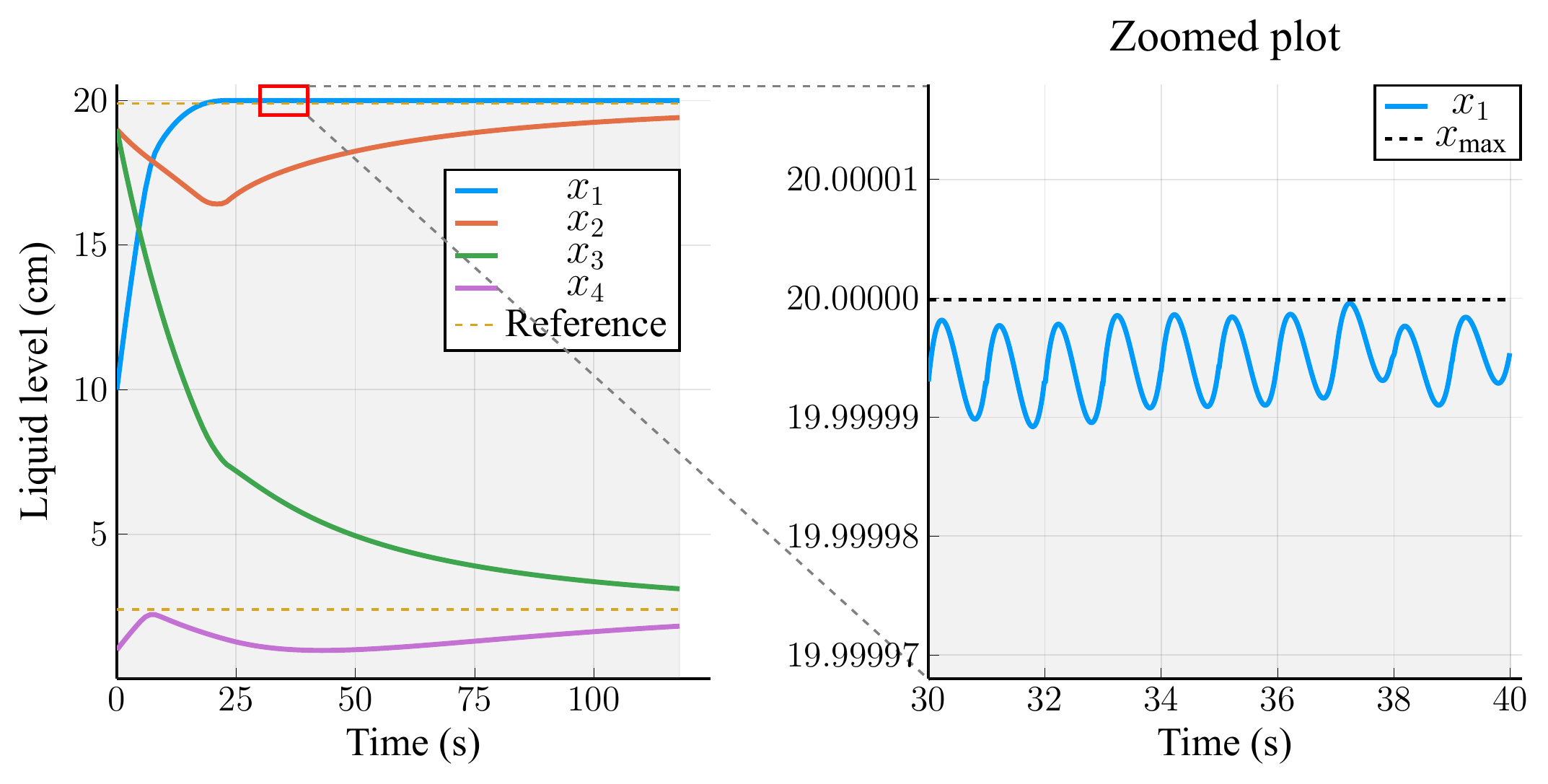}
	\caption{The states of continuous time linear MPC using our proposed solver. The gray area in the figures denotes the feasible region of states. The states using our proposed MPC input stays in the feasible region for whole time interval.}
	\vspace{-7pt}\end{figure}

As shown in Fig.~\ref{fig:state_input}, at the first glance, the state trajectories obtained from both solvers are nearly identical. However, upon close inspection, it can be seen that even though the QP-based controller satisfies the constraints at discrete-time instants, the constraints are violated in between sampling instants. In contrast, the proposed algorithm ensures constraint satisfaction on the whole time interval. 

\subsection{Computational Speed Performance}

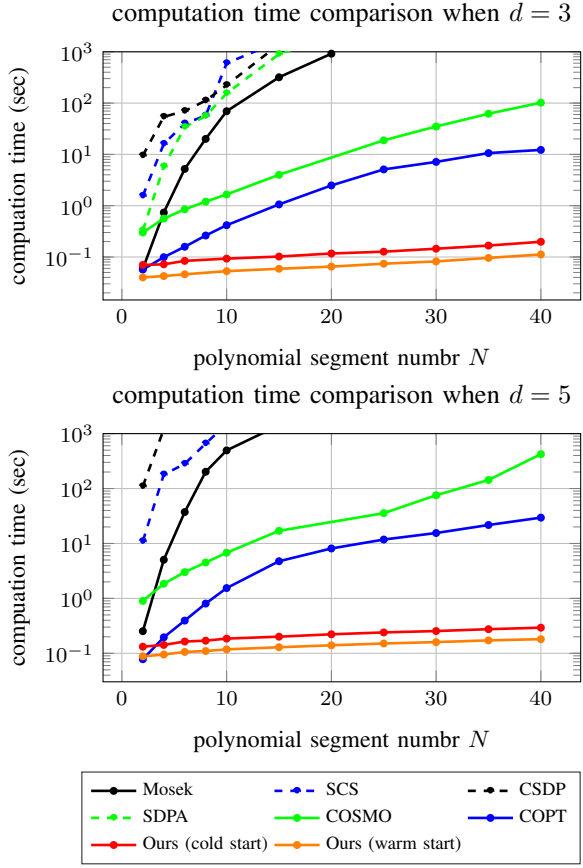
\begin{figure}[!htbp]
	\centering
	\begin{tikzpicture}
\begin{semilogyaxis}
	[width=2.5in,
	height=1.3in,
	at={(0in,2.0in)},
	scale only axis,
	ymin=0.0,
	ymax=1000,
	title={computation time comparison when $d=3$},
	xlabel={\small polynomial segment numbr $N$},
	ylabel={\small compuation time (sec)},
	y label style={at={(axis description cs:0.06,.5)},anchor=south},
	xticklabel style = {font=\footnotesize},
	yticklabel style = {font=\footnotesize},
	ytick={0.01,0.1,1,10,100,1000},
	axis background/.style={fill=white},
	xmajorgrids,
	ymajorgrids,
	legend style={at={(1.0,-1.9)}, anchor=north east, legend cell align=left, align=left, draw=white!15!black, legend columns=3, font=\scriptsize}]
		
	\addplot [color={black}, line width=1pt, mark=*, mark size=1pt]
	table[row sep={\\}]
	{   2  0.05896880999716814\\
		4  0.7353868739992322\\
		6  5.217759219001891\\
		8  20.04601449999973\\
		10 69.39\\
		15 318.5\\
		20 919.09\\
	}
	;\addlegendentry{Mosek}
	
		\addplot [color={blue}, dashed, line width=1pt, mark=*, mark size=1pt]
	table[row sep={\\}]
	{   
		2  1.56\\
		4  16 \\
		6  40\\
		8  57\\
		10 598 \\
		15 1596\\
	}
	;\addlegendentry{SCS}
	
		\addplot [color={black}, dashed, line width=1pt, mark=*, mark size=1pt]
table[row sep={\\}]
{   
	2  9.5\\
	4  54 \\
	6  70\\
	8  112\\
	10 224 \\
	15 1400\\
	20 4000\\
}
;\addlegendentry{CSDP}

		\addplot [color={green}, dashed, line width=1pt, mark=*, mark size=1pt]
table[row sep={\\}]
{   
	2  0.33\\
	4  5.8 \\
	6  34\\
	8  56\\
	10 155 \\
	15 900\\
	20 2000\\
}
;\addlegendentry{SDPA}

		\addplot [color={green},  line width=1pt, mark=*, mark size=1pt]
table[row sep={\\}]
{   
	2  0.30\\
	4  0.56 \\
	6  0.85\\
	8  1.2\\
	10 1.66 \\
	15 4.0
	20 9.3\\
	25 18.8\\
	30 35\\
	35 62\\
	40 102\\
}
;\addlegendentry{COSMO}

	\addplot [color={blue}, line width=1pt, mark=*, mark size=1pt]
	table[row sep={\\}]
	{   
		2  0.056915044784545\\
		4  0.09916400909423828\\
		6  0.15840888023376465\\
		8  0.2623138427734375\\
		10 0.41673707962036133\\
		15 1.0596790313720703\\
		20 2.484106874465942\\
		25 5.095527935028076\\
		30 7.146082878112793\\
		35 10.610167026519775\\
		40 12.233412981033325\\
	}
	;\addlegendentry{COPT}

		\addplot [color={red}, line width=1pt, mark=*, mark size=1pt]
	table[row sep={\\}]
{   2  0.07\\
	4  0.072\\
	6  0.084
	8  0.089\\
	10 0.093\\
	15 0.102\\
	20 0.117\\
	25 0.127\\
	30 0.145\\
	35 0.167\\
	40 0.198 \\
	}
	;\addlegendentry{Ours (cold start)}

		\addplot [color={orange}, line width=1pt, mark=*, mark size=1pt]
table[row sep={\\}]
{   2  0.04\\
	4  0.0425\\
	6  0.046
	8  0.049\\
	10 0.053\\
	15 0.059\\
	20 0.065\\
	25 0.074\\
	30 0.082\\
	35 0.096\\
	40 0.112 \\
}
;\addlegendentry{Ours (warm start)}

\end{semilogyaxis}


\begin{semilogyaxis}
	[width=2.5in,
	height=1.3in,
	at={(0in,0in)},
	scale only axis,
	ymin=0.0,
	ymax=1000,
	title={computation time comparison when $d=5$},
	xlabel={\small polynomial segment numbr $N$},
	ylabel={\small compuation time (sec)},
	y label style={at={(axis description cs:0.06,.5)},anchor=south},
	xticklabel style = {font=\footnotesize},
	yticklabel style = {font=\footnotesize},
	ytick={0.01,0.1,1,10,100,1000},
	axis background/.style={fill=white},
	xmajorgrids,
	ymajorgrids,
	legend style={at={(1.0,-0.38)}, anchor=north east, legend cell align=left, align=left, draw=white!15!black, legend columns=3, font=\scriptsize}]
	
	\addplot [color={black}, line width=1pt, mark=*, mark size=1pt]
	table[row sep={\\}]
	{   2  0.253139192\\
		4  5.044812105003075\\
		6  37.376431239998055\\
		8  202.04601449999973\\
		10 495.14\\
		15 1418.5\\
		20 1919.09\\
	}
	;
	
	\addplot [color={blue}, dashed, line width=1pt, mark=*, mark size=1pt]
	table[row sep={\\}]
	{   
		2  11\\
		4  180 \\
		6  281\\
		8  660\\
		10 1490 \\
		15 2596\\
	}
	;
	
	\addplot [color={black}, dashed, line width=1pt, mark=*, mark size=1pt]
	table[row sep={\\}]
	{   
		2  110\\
		4  1200 \\
	};
%
	
		\addplot [color={green},  line width=1pt, mark=*, mark size=1pt]
	table[row sep={\\}]
	{   
		2  0.9\\
		4  1.85 \\
		6  3.01\\
		8  4.5\\
		10 6.78 \\
		15 17.0\\
		25 35.6\\
		30 75.5\\
		35 143\\
		40 423\\
	}
	;
	
	\addplot [color={blue}, line width=1pt, mark=*, mark size=1pt]
	table[row sep={\\}]
	{   
		2  0.07804203033447266\\
		4  0.19564008712768555\\
	    6  0.3945600986480713\\
		8  0.805131196975708\\
		10 1.5411279201507568\\
		15  4.748100996017456\\
		20  8.111303091049194\\
		25 11.787981033325195\\
		30  15.489894151687622\\
		35  21.693773984909058\\
		40  29.536084175109863\\
	}
	;

	\addplot [color={red}, line width=1pt, mark=*, mark size=1pt]
	table[row sep={\\}]
	{   2  0.132\\
		4  0.1422\\
		6  0.164\\
		8  0.170\\
		10 0.1854\\
		15 0.201\\
		20 0.222\\
		25 0.240\\
		30 0.2544\\
		35 0.275\\
		40 0.293 \\
	}
	;
	
	\addplot [color={orange}, line width=1pt, mark=*, mark size=1pt]
	table[row sep={\\}]
	{   2  0.088\\
		4  0.095\\
		6  0.106\\
		8  0.110\\
		10 0.118\\
		15 0.129\\
		20 0.140\\
		25 0.151\\
		30 0.160\\
		35 0.172\\
		40 0.181 \\
	}
	;

\end{semilogyaxis}

\end{tikzpicture}
	\caption{The states of discrete time linear MPC using the QP solver.}\label{fig:time}
\end{figure}

In the following, we compare the computational speed performance of our proposed algorithm and several off-the-shelf solvers (on Problem \ref{pb:conic}) under different numbers of polynomial degrees $d$ and polynomial segments $N$. 
The block number for GPU acceleration is set as 128. The number of threads on every block is $\lceil \frac{pN}{128}\rceil$. The computational time in Figure \ref{fig:time} is the average solving time of the first 100 apply steps. The step sizes are $\alpha=0.2,\beta=0.4$. The computation platform is the same as in Subsection \ref{subsec:GPU}, i.e., a desktop computer equipped with an AMD Ryzen Threadripper 3970X 32-Core Processor and an NVIDIA GeForce RTX 3080 GPU. The real number calculations on GPU are floating point number with hybrid precision 32-bit and 16-bit, which is computationally efficient and accurate enough for control applications. As for comparison, the other solvers are of default precision 64-bit. Thus, the time comparison may not be equal but represents our computation speed superiority to some extent.

As shown in Fig. \ref{fig:time}, our proposed algorithm has better scalability for large problems (especially lagre $N$), and has low computational time promising for real-time control applications. The warm-start technique introduced in Subsection \ref{subsec:GPU} can effectively reduce the computation time by reducing iterations. For off-the-shelf solvers, COSMO and COPT perform well on large-scale problems compared to other solvers. However, their computation is still slow and incompatible with real-time control scenarios.

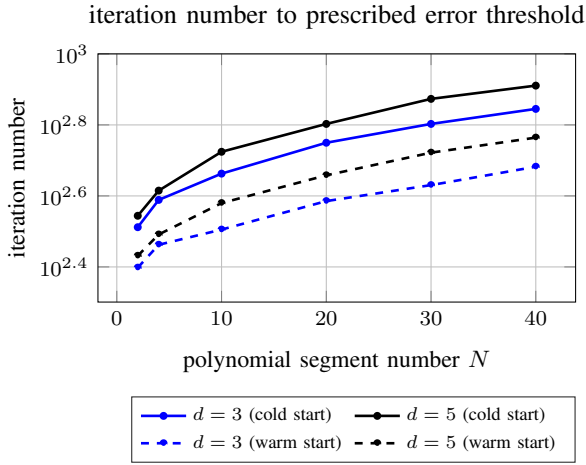
\begin{figure}[!htbp]
	\centering
	\begin{tikzpicture}

	\begin{semilogyaxis}
		[width=2.5in,
		height=1.3in,
		at={(0in,0in)},
		scale only axis,
		ymin=200,
		ymax=1000,
		title={iteration number to prescribed error threshold},
		xlabel={\small polynomial segment number $N$},
		ylabel={\small iteration number},
		y label style={at={(axis description cs:0.06,.5)},anchor=south},
		xticklabel style = {font=\footnotesize},
		yticklabel style = {font=\footnotesize},
		axis background/.style={fill=white},
		xmajorgrids,
		ymajorgrids,
		legend style={at={(1.0,-0.38)}, anchor=north east, legend columns=2, legend cell align=left, align=left, draw=white!15!black, legend columns=2, font=\scriptsize}]
		
		\addplot [color={blue}, line width=1pt, mark=*, mark size=1pt]
		table[row sep={\\}]
		{   2  325\\
			4  388\\
			10 460\\
			20 562\\
			30 635\\
			40 700 \\
		}
		;\addlegendentry{$d=3$ (cold start) }
		
		\addplot [color={black}, line width=1pt, mark=*, mark size=1pt]
		table[row sep={\\}]
		{   2  350\\
			4  412\\
			10 530\\
			20 635\\
			30 747\\
			40 814 \\
		}
		;\addlegendentry{$d=5$ (cold start) }

		\addplot [color={blue}, dashed, line width=1pt, mark=*, mark size=1pt]
		table[row sep={\\}]
		{   2  250\\
			4  290\\
			10 320\\
			20 385\\
			30 427\\
			40 481 \\
		}
		;\addlegendentry{$d=3$ (warm start) }
		
		\addplot [color={black}, dashed, line width=1pt, mark=*, mark size=1pt]
		table[row sep={\\}]
		{   2  270\\
			4  310\\
			10 380\\
			20 455\\
			30 527\\
			40 581 \\
		}
		;\addlegendentry{$d=5$ (warm start) }

		%
	\end{semilogyaxis}
	
\end{tikzpicture}
	\caption{The number of iterations required to reach residue $(\epsilon^k)_1<1\times 10^{-2}$, where $(\epsilon^k)_1$ means residue of segment 1. Warm start can effectively reduce the number of iterations.}\label{fig:iters}
\end{figure}
We demonstrate the number of iterations required to reach $(\epsilon^k)_1<10^{-2}$ for different problem sizes in Fig \ref{fig:iters}. The iteration number required grows gently as the problem size grows, which also corroborates the scalability of our proposed solver. The warm-start technique introduced in Subsection \ref{subsec:GPU} can effectively reduce the iteration number.

\begin{figure}[!htbp]
	\centering
	\input{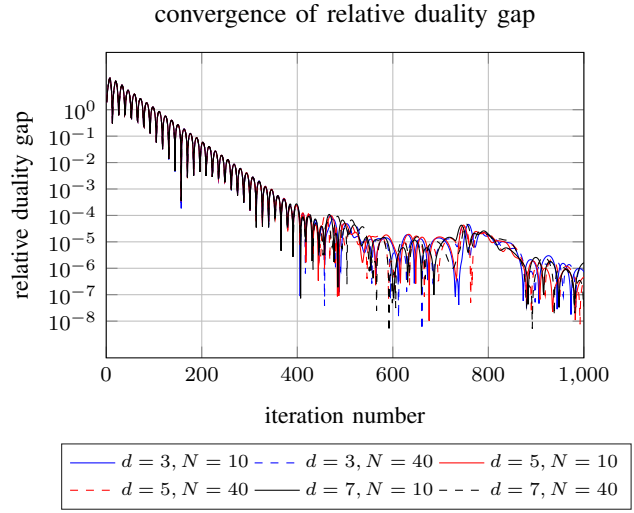}
	\caption{The relative duality gap of our proposed algorithm (cold start). Problem sizes scarcely influence the convergence speed of the relative duality gap. }\label{fig:converge}
\end{figure}
Define Lagrange function of Problem \ref{pb:conic} as $\mathscr{L}(\bc,s,\bX;\lambda)$, then the relative duality gap is defined as
$$\frac{1}{J(\bc)}\left[\inf_{\bc,s,\bX}\mathscr{L}(\bc,s,\bX;\lambda)-\sup_{\lambda}\mathscr{L}(\bc,s,\bX;\lambda)\right],
$$
where $J(c)$ is the objective value of Problem \ref{pb:poly_opt} and equivalently Problem \ref{pb:conic}.
We demonstrate the convergence of relative duality gap with respect to iteration number in Fig. \ref{fig:converge}. The relative duality gap converged below $10^{-5}$ within approximately 500 iterations. The problem sizes scarcely influence the convergence speed of the relative duality gap, which also indicates good scalability of our proposed algorithm.

		\section{conclusion}
		
		In this paper, we aim to address continuous-time path-constrained linear MPC problems while ensuring that path constraints are satisfied at every time interval. To achieve this, we propose an algorithm that utilizes differential flatness to eliminate dynamic constraints. Furthermore, by parameterizing the flat output with piecewise polynomials, we formulate a polynomial optimization problem where the decision variables are finite-dimensional polynomial coefficients, and the inequality path constraints are polynomial non-negativity constraints on intervals, which remain infinite-dimensional.
Taking advantage of the Markov-Luk\'{a}cs theorem from SOS theory, we transform the polynomial optimization problem into an equivalent SDP problem that is computationally tractable.
To accelerate the solving process of the SDP problem, we use a customized PDHG algorithm, which exploits the block-diagonal structure of the PSD matrix to perform paralleled computation. The numerical simulation of a quadruple-tank process validates that our proposed algorithm can ensure that the path constraints are satisfied at every time interval. Moreover, the parallel accelerated design of our algorithm results in superior computational speed performance.


		
		\bibliographystyle{IEEEtran}
		\bibliography{ref_cdc2023}
		
	\end{document}